# Effect of Encoding Method on the Distribution of Cardiac Arrhythmias

Luis A. Mora & Jhon E. Amaya

*Abstract*—This paper presents the evaluation of the effect of the method of ECG signal encoding, based on nonlinear characteristics such as information entropy and Lempel-Ziv complexity, on the distribution of cardiac arrhythmias. Initially proposed a procedure electrocardiographic gating to compensate for errors inherent in the process of filtering segments. For the evaluation of distributions and determine which of the different encoding methods produces greater separation between different kinds of arrhythmias studied (AFIB, AFL, SVTA, VT, Normal's), use a function based on the dispersion of the elements on the centroid of its class, the result being that the best encoding for the entire system is through the method of threshold value for a ternary code with E = 1 / 12.

## I. INTRODUCCION

In many cases to analyze ECG signals are used as measuring parameters the entropy and complexity of segments of such signals, but for the extraction of these features is necessary to convert the segments into sequences of symbols using different encoding methods, such as the method the threshold level [1-3] or by the slope of the signal [4-6], which can generate sequences composed of binary and ternary alphabets. Each method produces different sequences for the same segment of the signal, so that the extracted features of these sequences would have different values, however, the literature has not found a criterion for selecting the encoding method to use, as the type of analysis to be conducted.

Then analyzed the effect of these encoding methods on the distribution of different cardiac arrhythmias such as atrial flutter (AFL), supraventricular tachyarrhythmia (SVTA), ventricular tachycardia (VT) and atrial fibrillation (AFIB), which are extracted from ECG signals of the database of MIT-BIH [7]. This document is divided into the following sections. In the processing section details the electrocardiogram signal processing to extract the entropy and complexity. The following section describes the different encoding methods used, including key characteristics and relationships. In the evaluation section of the distribution, we present a method to estimate the quality of the distribution of the different arrhythmias using each method of coding. In the results section, presents the results obtained to evaluate the effect of coding on the distribution of arrhythmias, and finally presents the conclusions.

## II. SIGNAL PROCESSING

### A. Filtering

Before encoding signal segments for subsequent feature extraction, it is necessary to pre-process the signals to eliminate the effect of the various disturbances that affect the signal [8-9], such as breath noise, change at baseline, electromyographic noise, among others. This is done using the digital low-pass filters and high-pass, designed by Apaclia in [10], which are shown in (1) and (2) respectively.

$$F_{LP}(z) = \frac{1 - 2z^{-6} + z^{-12}}{36 - 72z^{-1} + 36z^{-2}} \quad (1)$$

$$F_{HP}(z) = \frac{-1 + 32z^{-16} - 32z^{-17} + z^{-32}}{32 - 32z^{-1}} \quad (2)$$

By combining these transfer functions are generated by the bandpass filter shown in (3), whose magnitude and phase response shown in Fig. 1.

$$F_{BP}(z) = \frac{(1 - 2z^{-6} + z^{-12})(-1 + 32z^{-16} - 32z^{-17} + z^{-32})}{1152(1 - 3z^{-1} + 3z^{-2} - z^{-3})} \quad (3)$$

However, to apply the filter described by equation (3), using the typical algorithm for digital filters, data loss occurs in each segment due to initial error and group delay [11-12].

$$y[n] = \sum_{i=0}^{M} b_i x[n-1] - \sum_{j=1}^{N} a_i y[n-j] \quad (4)$$

The initial error of a digital filter is that as shown in equation (4) the filter output signal *y[n]*, at any given time depends on previous values from itself and the input signal *x[n]*, values than at the beginning of the segment can not be known, so that the output of the filter is wrong in the first $l_i$ samples.

By another hand, the loss of information generated by group delay is due to the phase response of the filter produces a delay generated by the loss of the last segment $l_f$ samples. Mathematically group delay is defined as the derivative of the phase response of the filter on the frequency, as shown in (5), so if the filter is linear phase, then $l_f = \tau(\omega)$ and is constant for all ω.

---

·Manuscript received June 30, 2011.

　　Luis A, Mora is with the Laboratory of Instrumentation, Control and Automation, Universidad Nacional Experimental del Tachira, San Cristobal 5001, Venezuela (e-mail: lmmora@unet.edu.ve).

　　John E. Amaya is with the Laboratory of High Performance Computing, National University Experimental del Tachira, San Cristobal 5001, Venezuela (e-mail: jedgar@unet.edu.ve).

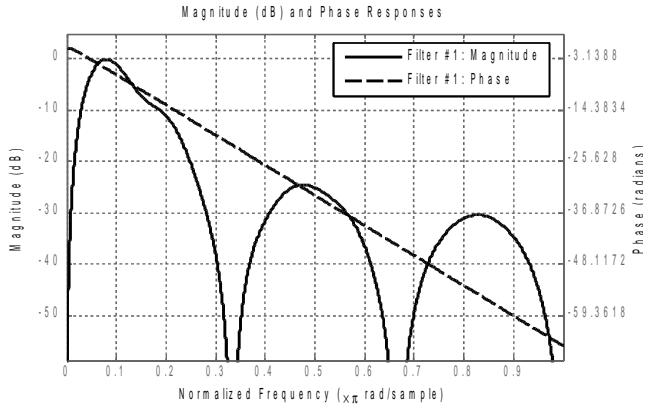

Fig. 1. Magnitude and phase responses of the bandpass filter used.

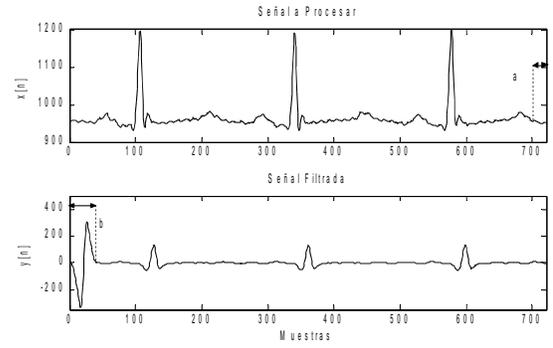

Fig. 2. Errors in the filter segment. (a) Information lost due to group delay. (b) Information loss due to initial error of the filter.

$$\tau(\omega)=\frac{-d(\sphericalangle H(j\omega))}{d\omega} \quad (5)$$

If the size of the segment to be processed is very large then these errors would be negligible, but for short length segments of information loss due to them is significant, especially if the filter order is high. In Figure 2 shows the effect of these errors to a segment of 721 samples of an ECG using the bandpass filter described in (3), which lost the last 21 samples due to the delay group and 44 samples due to initial error.

To compensate the errors described, we propose the following method: Let $H(z)$ a linear phase filter with group delay $l_f$ samples, where $M$ and $N$ represent the orders of the polynomials of the numerator and denominator of the filter transfer function, where $l_i = \max\{M, N\}$ and $X=\{X_1, X_2,..., X_n\}$ the segment of the signal to be processed. To correct the errors can follow the steps shown in Figure 3, which added to the original segment, the segment's $I=\{I_1, I_2,..., I_{ki}\}$ and $F=\{F_1, F_2,..., F_{kf}\}$ to compensate the initial error and group delay respectively, where $\forall I = X_1$ and $\forall F = X_n$, being $k_i \geq l_i$ and $k_f \geq l_f$. Then this new segment is filtered using the formula presented in (4) and later extract the segment with all the desired information. As shown in Figure 4, this method effectively eliminates the initial error and makes the filtered signal in the image below, is in phase with the original signal, thus compensating the group delay, the only disadvantage of this methodology to the processing time is greater than the normal filtering process, given the fact the original add additional segments and process a larger amount of data.

*B. Features:*

Usually for the analysis of electrocardiographic signals (ECG) are extracted temporal and morphological features of the electrical signals of the heart [13-15], being necessary for the identification of abnormalities usually 3 or more of these characteristics. However, another frequently used method is the extraction of nonlinear features such as entropy and complexity for carrying out the analysis [4], [5], [16-18]. In this paper, the extracted features are the Shannon entropy and complexity of Lempel-Ziv, which are used to evaluate the distribution of different electrocardiographic signals according to the encoding methods.

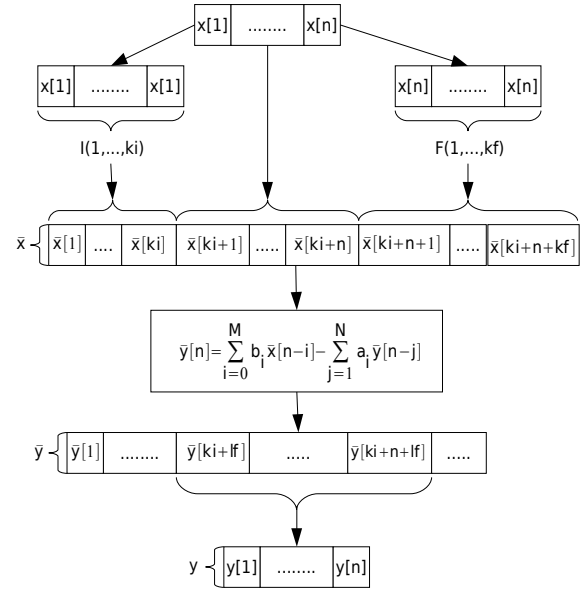

Fig. 3. Method used to compensate for the errors generated by digital filtering of a signal.

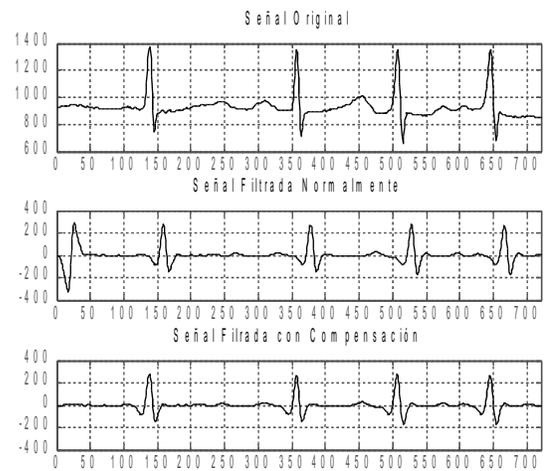

Fig. 4. Comparison of the proposed filtering method with respect to common filtering procedure.

*Shannon Entropy*

Also known as information entropy, the Shannon entropy is a mathematically rigorous measure of the lack of information associated with a probability distribution of a data set [19-21]. Was defined by Shannon in [22] as shown in equation (6), where α is the size of the alphabet of possible values that can take a data set *S* and $p_i$ the probability that a random sample of *S* has the value $S_i$. Whose normalized form is given by (7).

$$H = -\sum_{i=0}^{\alpha} p_i \log_2(p_i) \qquad (6)$$

$$H_n = \frac{-\sum_{i=0}^{\alpha} p_i \log_2(p_i)}{\log_2(\alpha)} \qquad (7)$$

*Complexity of a Sequence:*

There are several techniques to measure the complexity of a sequence, one of them is the complexity of Lempel-Ziv or LZ complexity, which is frequently used in cryptography and data compression [23]. This measure of complexity "scans" a sequence *S* of length *n* from left to right and added a new word each time a memory is a sub-sequence of consecutive digits has not been found previously, using complexity metrics to measure the size of the alphabet and the number of words found in *S* [24]. That is, the complexity *c(s)* of a sequence *S* = {$S_1, S_2,..., S_n$} consists of elements of an alphabet *A* and the minimum number of steps or sub-sequences necessary to construct *S*.

There are several algorithms to compute the Lempel-Ziv complexity, this study used a variation of the algorithms presented in [3], [25-26]. But a comparison of complexity between different sequences using the normalized complexity [1], [3], [25].

$$C(s) = \frac{c(s)}{b(s)}; b(s) = \frac{n}{\log_\alpha(n)} \qquad (8)$$

*C. Selecting the length of the segment.*

For this study, the ECG signals are processed by segment, so that a proper selection of the length of it is of vital importance, as being very short features would not be very representative and ultimately the cost of processing would be very high. To select the minimum length is taken into account the complexity according to [24] as,

$$c(s) = \frac{n}{(1-\varepsilon_n)\log_\alpha(n)}; \varepsilon_n = 2\frac{1+\log_\alpha(\log_\alpha(\alpha n))}{\log_\alpha(n)} \qquad (9)$$

$$\frac{1+\log_\alpha(\log_\alpha(\alpha n))}{\log_\alpha(n)} < \frac{1}{2} \qquad (10)$$

This indicates that according to (9) for the measurement of c(s) must be met to make sense $\varepsilon_n < 1$, being that $\varepsilon_n$ depends on the size of the alphabet and the length of the segment, then the complexity measure is valid if equation is satisfied (10).

So the length of a segment composed of a binary alphabet is $n \geq 361$ and for a ternary alphabet $n \geq 366$. In addition to study arrhythmias are sustained arrhythmias with frequencies above 100 beats per second, and in the case of ventricular tachycardia (VT) arrhythmia must be 3 beats or more [27]. Whereas electrocardigráficas signals found in MIT-BIH [7] have been acquired at a frequency of 360 samples per second, the segments must be longer than 648 samples to cover the three keystrokes required to identify a VT with a minimum frequency of 100 beats / sec. According to this study were selected segments of two seconds duration, ie 720 samples to meet both constraints.

### III. Encoding Methods

To extract features is first necessary to encode the segment *Y* of the filtered signal of the electrocardiogram, for it uses the techniques of the slope and the threshold value, which can generate binary or ternary codes. In the case of the slope method the encoding is based on the assignment of a symbol of the alphabet used according to the difference between two consecutive samples of the signal, which can be positive, negative or zero, and has the advantage reduces the length of the segment in a sample, but is very sensitive to noise. The threshold value methodology codes allocation is done by comparing the signal samples over a certain value, usually determined by the average value of the signal. In this work, to build the binary alphabets used the symbols "0" and "1", and for ternary symbols and "-1", "0" and "1". Table 1 presents a summary of the coding algorithms used to assign these symbols according to the methods proposed. As *S* the coded signal, *P(k)* the slope of the signal between samples *k* and *k*+1, *T* the threshold value, *Y* the mean value of the segment and *E* deviation from the mean maximum and minimum value of the signal. An example of the effect of these methods is shown in fig. 4 and 5, where an electrocardiogram of a normal rate of 721 samples length is encoded using the slope for binary and ternary alphabet (Fig. 4b and 4c respectively) and using the threshold value for codes 2 and 3 symbols (5b and 5d, respectively). In the latter method can be seen in Fig. 5a and 5c the reference values used in each case, a set of line segments. You can also distinguish the differences in the shapes of each coded segment, since each method produces different sequences.

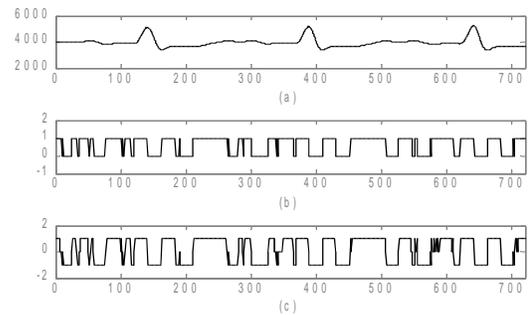

Fig. 5. Example of coding using the slope method (a) Normal ECG, (b) binary encoding, (c) ternary coding

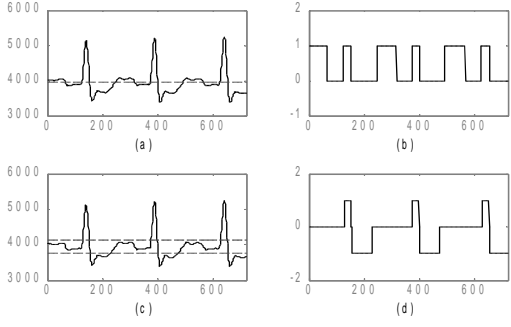

Fig. 6. Example of coding using the threshold method, (a) and (c) normal ECG, (b) binary encoding, and (d) ternary coding

TABLE 1: ENCODING METHODS

| Method | Principle | Encoding |
|---|---|---|
| Slope Binary | $P(k)=\Delta y(k+1,k)$ | $S(k)=1 ; \forall k \in \mathbb{N}, P(k) \geq 0$ <br> $S(k)=0 ; \forall k \in \mathbb{N}, P(k) < 0$ |
| Slope Ternay | $P(k)=\Delta y(k+1,k)$ | $S(k)=1 ; \forall k \in \mathbb{N}, P(k) > 0$ <br> $S(k)=0 ; \forall k \in \mathbb{N}, P(k) = 0$ <br> $S(k)=-1 ; \forall k \in \mathbb{N}, P(k) < 0$ |
| Threshold Binary | $T=Y+E(máx(Y)-mín(Y))$ | $S(k)=1 ; \forall k \in \mathbb{N}, y(k) \geq T$ <br> $S(k)=0 ; \forall k \in \mathbb{N}, y(k) < T$ |
| Threshold Ternary | $Ta=Y+E(máx(Y)-mín(Y))$ <br> $Tb=Y-E(máx(Y)-mín(Y))$ | $S(k)=1 ; \forall k \in \mathbb{N}, y(k) > Ta$ <br> $S(k)=0 ; \forall k \in \mathbb{N}, Tb \leq y(k) \leq Ta$ <br> $S(k)=-1 ; \forall k \in \mathbb{N}, y(k) < Tb$ |

## IV. EVALUATION OF THE DISTRIBUTION

To evaluate which of the encoding methods shown in Table 1 generates a better distribution of the features, then proposed an assessment based on the centroids of each class of data and the distance of each data to these centroids. To this end consider the following: Let $D$ be a set of data distributed in $m$ class, $D=\{C_1, C_2,..., C_m\}$, where each class $C_i$ contains $n_i$ elements $C_i=\{e_{i1}, e_{i2},..., e_{ini}\}$, and there a centroid $\Phi_i$ for each class.

Lemma 1: Can be evaluate the distribution of elements in the data set D from the ratio of the distance between each element to the centroid of its own class, with respect to its distance from the centroid of the other classes. Considering a spatial centroid as the center of the distribution of the elements of a class, its centroid is determined by:

$$\Phi_i = \frac{\sum_{j=1}^{n_i} e_{ij}}{n_i} \quad (11)$$

*Definition 1*: Let $l_i$ the distance between an element $e_{ij}$ and the centroid $\Phi_i$, the distribution function $\lambda$ of a class of elements of a data set is defined as the number of elements of the class whose centroid distance is greater or equal to the distance from the centroid of any other class. That is:

$$\lambda_i = \sum_{j=1}^{n_i} \xi(e_{ij}) ; \xi(e_{ij}) = \begin{cases} 1, & l_i \geq l_h, \forall h \neq i \\ 0, & otherwise \end{cases} \quad (12)$$

The function $\lambda_i$ is annulled when all the elements of the class Ci are clustered near its centroid, and is maximal when there is superposition or total overlap between 2 or more classes. The maximum value that can take the distribution function of a class is equal to the number of elements of that class, so the normalized distribution of a class is given by (13).

$$\lambda_{in} = \frac{\lambda_i}{n_i} \quad (13)$$

*Definition 2*: It is understood as the distribution of a data set $\lambda_D$, the sum of the distributions of each class is made up, as shown in (14).

$$\lambda_D = \sum_{i=1}^{m} \lambda_i \quad (14)$$

*Definition 3:* this refers to a normalized distribution $\lambda_{DN}$ the relationship between the distribution system and the total number of elements that comprise it.

$$\lambda_{DN} = \frac{\lambda_D}{N} ; N = \sum_{i=1}^{m} n_i \quad (15)$$

*Definition 4*: is understood as the average distribution of a data set $\lambda_{DP}$, the relationship between the sum of the normalized distributions of the classes that make up the system and the total number classes.

$$\lambda_{DP} = \frac{\sum_{i=1}^{m} \lambda_{in}}{m} \quad (16)$$

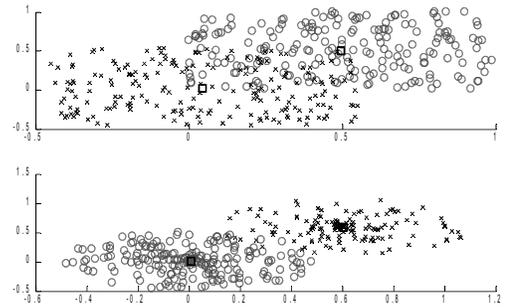

Fig. 7. Two examples of data distributions

In Figure 7 are two examples of distributions to which they apply this method of evaluation. At the top of the data set consists of two classes with 200 items each and are slightly overlapping in the image below the data set also has two classes, but one has 200 elements and the other only 100, having a greater apart from each other. Table 2 gather the results of evaluating these distributions, note that the normalized distributions of the classes of the top image in Figure 7 are higher than those in the bottom, this is because the first case classes are more superimposed. Similarly $\lambda_D$ is higher due to the separation between classes. Additionally $\lambda_{DN}$ and $\lambda_{DP}$ are equivalent when all classes have the same number of elements.

TABLE 2: EVALUATION OF THE DISTRIBUTIONS PRESENTED IN FIG. 7

| Dist. | $\lambda_1$ | $\lambda_{1_n}$ | $\lambda_2$ | $\lambda_{2_n}$ | $\lambda_D$ | $\lambda_{DN}$ | $\lambda_{DP}$ |
|---|---|---|---|---|---|---|---|
| top | 27 | 0,135 | 26 | 0,130 | 53 | 0,133 | 0,133 |
| lower | 3 | 0,015 | 4 | 0,027 | 7 | 0,018 | 0,021 |

## V. RESULTS

For the study decided to work with a deviation from the average value of E={-1/10, -1/20, 1/20,1/10} in case of binary encoding for the method of threshold value, and E={1/8, 1/10, 1/12, 1/14, 1/16, 1/20} in case of ternary encoding the same method. The selection of these values is that for values outside this range the encoder generates no significant variations in the extracted features.

To evaluate the effect of different encoding methods on the distribution of arrhythmias were extracted 2124 segments of the database electrocardiographic signals from MIT-BIH, distributed as follows: 800 Normal, 325 AFIB, AFL 324, 325 SVTA and 350 VT.

In applying the definitions 3 and 4 of Section III, we found that the encoding produces better distribution for the entire system is one that uses the ternary threshold value E = 1 / 12 with λDN = 0, 4035 and λDP = 0.4358.Table 3 can be best encoding methods as evaluated using the definition 2 of Section III. For more details on the evaluation of the distributions generated by the different methods you can visit the web address presented by Mora and Amaya [28].

TABLE 3: ENCODING METHODS WITH BETTER DISTRIBUTION FOR EACH PAIR OF ARRHYTHMIAS

| Pairs of arrhythmias | Method | $\lambda_{in}$ |
|---|---|---|
| Normal Vs Afib | PT[1] | 0.1929 |
| Normal Vs Afl | VUB[2] con E=-1/20 | 0.1032 |
| Normal Vs Svta | VUB[2] con E=1/10 | 0.1067 |
| Normal Vs Vt | VUB[2] con E=1/10 | 0.1383 |
| Afib Vs Afl | VUT[3] con E=1/10 | 0.1156 |
| Afib Vs Svta | VUT[3] con E=1/8 | 0.2261 |
| Afib Vs Vt | VUB[2] con E=1/20 | 0.1896 |
| Afl Vs Svta | VUT[3] con E=1/12 | 0.0447 |
| Afl Vs Vt | VUB[2] con E=1/20 | 0.1231 |
| Svta Vs Vt | VUT[3] con E=1/10 | 0.2444 |

[1] Slope method with binary encoding
[2] Threshold method with binary encoding
[3] Threshold method with ternary coding

## VI. CONCLUSIONS

The function λ relates the distance from one element to the centroid of its class, with the distance to the centroid of another class, which allows an estimate of the separation between classes, ie the lower the value of λ, the higher the separation between classes, as seen in Figure 7, the lower classes have a separation greater than the upper classes, so that in evaluating their results are lower distributions as seen in Table 2 . Using the function λ can observe the effect of each encoding method on the distribution of arrhythmias in entropy-complexity plane, with the VUB and * STV methods, which create greater separation between classes and the least recommended is the encoding binary using the slope of the signal.

By selecting the correct encoding method can make a separation between the different types of arrhythmias that can discriminate each type using only the Shannon entropy and complexity of Lempel-Ziv, decreasing the amount of features extracted by Mahmoodabadi et al. in [15] for the classification of the signals.